\documentclass{article}
\usepackage{amsfonts,amsmath,amssymb,amsthm,cite,enumerate,graphicx,latexsym,mathrsfs}
\numberwithin{equation}{section}

\newtheorem{definition}{\indent Definition}[section]
\newtheorem{lemma}{\indent Lemma}[section]

\newtheorem{remark}{\indent Remark}[section]

\newtheorem{theorem}{\indent Theorem}[section]
\def\disp{\displaystyle}
\def\a{\alpha}

\def\r{\rho}

\def\ld{\lambda}

\def\l{\leq}
\def\g{\geq}

\def\i{\infty}

\def \B{\begin{equation}}
\def\bbb{\begin{array}{lll}}

\def\eee{\end{array}}
\def \E{\end{equation}}

\def\it{\indent}

\title{{Singularity formation for compressible Euler equations with time-dependent damping
\thanks{This work is supported in part by the National Natural Science Foundation of China (Grant No. 11671237)}}}

\author{{Ying Sui\thanks{Email: suiying4320@163.com}}~~~~{Huimin Yu\thanks{Corresponding author, Email: 106139@sdnu.edu.cn}}\\
\small\textit{$$ School of Mathematics and Statistics, Shandong Normal University, Jinan 250014, China.}\\}

\begin{document}
\date{}
\maketitle

\textbf{Abstract:}
In this paper, we consider the compressible Euler equations with time-dependent damping $\frac{\a}{(1+t)^\lambda}u$ in one space dimension. By constructing  ¡°decoupled¡± Riccati type equations for smooth solutions,
we provide some sufficient conditions under which the classical solutions must break down in finite time.
As a byproduct, we show that the derivatives blow up, somewhat like the formation of shock wave,
if the derivatives of initial data are appropriately large at a point even
when the damping coefficient goes to infinity with a algebraic growth rate.
We study the case $\lambda\neq1$ and $\lambda=1$ respectively,
moreover, our results have no restrictions on the size of solutions and the positivity/monotonicity of the initial Riemann invariants.
In addition, for $1<\gamma<3$ we provide time-dependent lower bounds on density for arbitrary classical solutions,
without any additional assumptions on the initial data.

\textbf{Keywords:} Singularity formation, compressible Euler equations, time-dependent damping, shock wave.

\section{Introduction}

\ \ \ \
In this paper, we consider the one dimensional compressible
Euler equations with time-dependent damping in Lagrangian coordinates:
\begin{equation}\label{1}
\left\{
\begin{array}{l}
\tau_t-u_x = 0,\\
u_t+ p_x = -\frac{\a}{(1+t)^\lambda}u,
 \end{array}
 \right.
\end{equation}
where $\tau=\tau(x,t)$ and $u=u(x,t)$ are the specific volume and velocity of the flow at location $x\in R$ and time $t\in R^+$. For simplicity,
we assume the gas is ideal polytropic and the gas pressure
\begin{equation}\label{3}
\begin{aligned}
p = K\tau^{-\gamma} \ \ \text{for constants $K>0$ and} \ \ \gamma>1.
 \end{aligned}
\end{equation}
Besides, $\a\geq 0,\ \lambda\in R$ are two constants,
and the term $-\frac{\a}{(1+t)^\lambda}u$ is the so-called damping effect on the fluid.
Physically, model (\ref{1}) is used to describe the compressible isentropic flow through porous medium with unsteady drag force. In this paper, we want to investigate the blow up phenomena of system (\ref{1}) with the $C^1$ initial data
\begin{equation}\label{4}
\begin{aligned}
(\tau, u)(x, 0) = (\tau_0, u_0)(x).
 \end{aligned}
\end{equation}
In other words, we are concerned with the conditions (on the initial data) under which a shock forms in a
classical $C^1$ solution.

When $\a=0$, the system $(\ref{1})$ reduces to the standard compressible Euler equations, which has been studied extensively.
The smooth solutions blow up in general due to the formation of shocks (\cite{Dafermos,Lax,Smoller,Slemrod}).
Especially, \cite{Chen3} give a sufficient and necessary condition for the formation of blow up phenomena.
During the proof of \cite{Chen3}, the uniform lower bound estimates of density plays an important role. \\
\it When $\alpha>0,$  $ \lambda=0$, the system $(\ref{1})$ turns into the compressible Euler equations with constant coefficient damping.
In 1988, Lin\cite{Lin} considered the singularity formation mechanism and the global existence of smooth solution for big initial data.
However, in Lin's proof, an additional condition are assumed on the initial data to ensure that the $C^1$-solution is strictly away from vacuum.
 There are also other works on the relaxation limit and asymptotic behavior in 1-D or even 3-D case, we refer to \cite{Hsiao,Mar1,Mar2,Mar3,Nish1,Nishi2,Sid} and the references therein.  \\
\it As for the time-dependent damping model (\ref{1}), Pan \cite{Pan1, Pan2} gave the thresholds of $\a$ and $\ld$
to separate the existence and nonexistence of global smooth solutions in small data regime.
Sugiyama \cite{Sugi} obtained the sharp upper and lower estimates of life span for some cases of $\ld$ and $\a$.
Recently, \cite{Chens} gave some sufficient conditions to make all solutions blow up with monotonic initial Riemann invariants.
In particular, the blow up phenomena can be seen for small initial data. We can also refer \cite{Hou1,Hou2,Cui} for more interesting topics (such as the existence of smooth solutions and their approximate behavior) on this model.\\
\it About the lower bound estimates of density,\cite{Chens} assume the initial Riemann invariants
are bigger than some positive constant $\varepsilon_0$, then the density is away from vacuum, see Theorem 1.2 in \cite{Chens} for detail.
A straight forward question is: How about the case for general bounded initial Riemann invariants?
It is well known that there is no uniform lower bound of density for the compressible Euler equations with $L^\i$ initial data.
For example, the density may goes to zero with some time decay rate, see \cite{Chen2}.\\
\it In this paper, for the Euler equations with time-dependent damping,
we will give the lower bound estimates of density when $1<\gamma <3$, on which the blow up mechanism is discussed.
We also consider the singularity formation for $\gamma >3$.
We shall show that the derivatives of the smooth solutions to problem $(\ref{1})-(\ref{4})$ blow up if the derivatives of
initial data are appropriately large at a point even when the damping coefficient goes to infinity with algebraic growth rate,
which means the increasing damping effect can not cancel the hyperbolic effect of the Euler equations totally.
It is worthwhile pointing out that there have no any monotonic assumption on the initial Riemann invariants.
In the proof, we first consider the system (\ref{1}) in the form of Riemann invariants and differentiate it by space variables $x$,
then make a series of nonlinear changes of variables to get an uncoupled pair of Riccati type ODEs along characteristics.
Next, we  investigate the decoupled ODEs  to gain some sufficient conditions for singularity formation of the system (\ref{1}) in finite time.
Especially, to exhibit the blow up mechanism we need to obtain upper and lower bounds for density, which is a crucial key in the proof. Fortunately, inspired by the work of
\cite{Chen2}, we give the specific form of upper and lower bounds for smooth density.

The paper is organized as follows: in Section 2, we obtain a priori $L^\infty$
bounds for any weak solutions (including smooth solutions of course) to (\ref{1})-(\ref{4}).
We prove this result by borrow the method  established in \cite{Huang} for constant damping Euler equations.
In section 3, we consider the singularity formation for compressible Euler equations when $\ld\neq 1.$
While the blow up mechanism for $\ld=1$ are investigated in Section 4. In Section 5, we give some comments and discuss some further problems.

\section{Invariant Region and some preparation}
\ \ \ \
In this section, we will give an invariant region theorem for all weak solutions of
compressible Euler equations with time dependent-damping.
 The proof is very similar with the system with constant coefficient damping,
 which was considered in \cite{Huang}. Here we only give the sketch of the analysis. \\
\it Consider problem (\ref{1})-(\ref{4}) under the Euler coordinates, that is
\begin{equation}\label{E1}
\left\{
\begin{array}{l}
\r_t+(\r u)_x = 0,\\
(\r u)_t+ (\r u^2+\textsc{p}(\r))_x = -\frac{\a}{(1+t)^\lambda}\r u,
 \end{array}
 \right.
\end{equation}
with the initial data
\begin{equation}\label{E1-I}
\r(x,0)=\r_0,~ u(x,0)=u_0,
\end{equation}
where $\r=v^{-1},~\textsc{p}(\r)=K \r^\gamma$.
Borrow the method introduced in \cite{Huang}\cite{Fang},
we give an invariant region theorem for any $L^\i$ weak solutions to system (\ref{E1}).\\
\it To proceed the analysis, we first give the definition of $L^\i$ weak entropy solutions to system (\ref{E1}) and recall some fundamental results on the entropies for the standard Euler equations:
 \begin{equation}\label{E2}
\left\{
\begin{array}{l}
\r_t+(\r u)_x = 0,\\
(\r u)_t+ (\r u^2+\textsc{p}(\r))_x =0.
 \end{array}
 \right.
\end{equation}
\begin{definition} Denote $m=\r u$,
for any $T > 0$, the bounded measurable functions $(\r, m)(x,t) \in L^\i(R\times[0,T])$ are called weak entropy solutions of (\ref{E1})-(\ref{E1-I}), if
(\ref{E1})(\ref{E1-I}) and \begin{equation}\label{E1-E}
\eta_{t}+q_{x}+{\a \over {(1+t)^\lambda}}\eta_{m}m\leq0
\end{equation}
hold in the sense of distributions, where $(\eta,q)$ is any weak convex entropy-flux
pair $(\eta(\r,m), q(\r,m))$ satisfying
$\nabla q = \nabla \eta\nabla f,~ f=(m, {{m^2}\over {\r}}+K\rho^{\gamma})^{T},$ and
$\eta(0,0) = 0.$
\end{definition}
 As noted in \cite{Lions1}, all weak entropy-entropy flux pairs of $(\ref{E2})$ can be demonstrated by the following formulas:
\begin{equation}\label{*.1}
\begin{array}{l}
\eta(\rho,u)=\disp\int g(\xi)\chi(\xi;\rho,u)d\xi=\rho\int^{1}_{-1}g(u+z\rho^{\theta})(1-z^{2})^{\lambda}dz,\\
q(\rho,u)=\disp\int g(\xi)(\theta\xi+(1-\theta)u)\chi(\xi;\rho,u)d\xi
         =\rho\int^{1}_{-1}g(u+z\rho^{\theta})(u+\theta z\rho^{\theta})(1-z^{2})^{\lambda}dz,
\end{array}
\E
where $\theta=\disp\frac{\gamma-1}{2}$, $\lambda=\disp\frac{3-\gamma}{2(\gamma-1)}$, $g(\xi)$ is any smooth function of $\xi$,
and
\B\label{*.2}
\chi(\xi;\rho,u)=(\rho^{\gamma-1}-(\xi-u)^{2})^{\lambda}_{+}.
\E
The formula $(\ref{*.1})$ can be deduced from the entropy equation $\nabla q=\nabla\eta \nabla f$
by exploring the fundamental solutions of linear wave equations or the kinetic formulation. We remark that when $g(\xi)={1\over 2}\xi^2$, then $\eta={1\over 2}\r u^2+{K\over {\gamma-1}}\r^\gamma$ is mechanical energy.
In this section, like \cite{Huang}, we choose $g(\xi)=g_{k}(\xi)=e^{k\xi^{2}}$ in $(\ref{*.1})$,
then the corresponding entropy-entropy flux pair is
\begin{equation}\label{*.3}
\begin{array}{l}
\eta_k=\disp\rho\int^{1}_{-1}e^{k(u+z\rho^{\theta})^{2}}(1-z^{2})^{\lambda}dz,\\
q_k=\disp\rho\int^{1}_{-1}e^{k(u+z\rho^{\theta})^{2}}(u+\theta z\rho^{\theta})(1-z^{2})^{\lambda}dz.
\end{array}
\end{equation}
Because of the results given by P. L. Lions \textit{et~al.}\cite{Lions1},
the mechanical energy $\eta_{k}>0$ is convex obviously.
Therefore, for any $L^{\infty}$ weak entropy solution $(\rho,m)$, the entropy inequality
\begin{equation}\label{*.4}
\eta_{kt}+q_{kx}+{\a \over {(1+t)^\lambda}}\eta_{km}m\leq0
\end{equation}
satisfies in the sense of distributions.
Using the same calculations in \cite{Huang}, we get the uniform estimate
 \begin{equation}\label{E3}
 \begin{array}{l}
 \disp\max\{\|(u-\r^\theta)(x,t),(u+\r^\theta)(x,t)\|_{L^\i}\}\\\\
 \l \disp\max\{\|(u-\r^\theta)(x,0),(u+\r^\theta)(x,0)\|_{L^\i}\}
 \end{array}
 \end{equation}
for any $t>0$, providing the condition
\begin{equation}\label{*.5}
{\a \over {(1+t)^\lambda}}\eta_{km}m\geq 0
\end{equation}
meets. While the inequality $\eta_{km}m\geq 0$ have been checked in \cite{Huang}.
 Noticing the positivty of ${\a \over {(1+t)^\lambda}}$, we conclude the following Lemma:
\begin{lemma}\label{L2.1}
Suppose that $(\r_0, u_0)(x) \in  L^\i(R)$ satisfies
$$0 \leq \r_0(x) \leq C_0, \ \ \ \ |m_0(x)| \leq C_0\r_0(x), $$
for some positive constant $C_0>0$.
Let $(\r,u) \in L^\i(R\times [0,T])$ be any $L^\i$ weak entropy solution of the system $(\ref{E1})$ with
$ \gamma>1$. Then $(\r,m)$ satisfies
\begin{equation}\label{E4}
 0\l \r(x,t)\l \tilde{C}_0,~~ |m(x,t)|\l \tilde{C}_0\r(x,t),
\end{equation}
where the constant $\tilde{C}_0$, which can be choosed bigger than $\max\{C_0+C_0^\theta, (C_0+C_0^\theta)^{1\over \theta}\}$, depends solely on the initial data.
 \end{lemma}
Base on the above Lemma, the uniform bounded estimates for smooth solutions to system (\ref{1})-(\ref{4}) can be gained directly.
\begin{theorem}\label{t2.1}
Let $(v_0, u_0)(x)$ satisfies
\begin{equation} \label{E5*}
C_0^{-1} \leq v_0(x), ~~|u_0(x)| \leq C_0,
\end{equation}
and $(v,u) \in C^1(R\times [0,T])$ be the smooth solution to system $(\ref{1})-(\ref{4})$ for $T>0$.
Then $(v,u)(x,t)$ is uniformly bounded in the following form:
\begin{equation}\label{E5}
\tilde{C}_0^{-1}\leq v(x,t),~~ |u(x,t)|\l \tilde{C}_0,
\end{equation}
where $\tilde{C}_0$ is the same as in Lemma \ref{L2.1}.
 \end{theorem}
\begin{remark}
Compaired with the result of Theorem 1.2 in \cite{Chens}, we generlize the assumptions on the initial Riemann invariants by (\ref{E5*}).
However, the uniform upper bound estimates of specific volume, i.e. the lower bound of density is absent.
We will use anther method to get the estimates in Section 3.2.
\end{remark}

Come back to the Lagrangian coordinate, denote
\begin{equation}\label{10}
\begin{aligned}
\phi :=\int^\infty_\tau c d\tau=\frac{2\sqrt{K\gamma}}{\gamma-1}\tau^{-\frac{\gamma-1}{2}}>0,
 \end{aligned}
\end{equation}
where the nonlinear Lagrangian sound speed $c$ is
\begin{equation*}\label{11}
\begin{aligned}
c :=\sqrt{-p_\tau} = \sqrt{K\gamma}\tau^{-\frac{\gamma+1}{2}}.
 \end{aligned}
\end{equation*}
It follows that
\begin{equation*}\label{12}
\begin{aligned}
\tau= K_\tau \phi^{-\frac{2}{\gamma-1}},
 \end{aligned}
\end{equation*}
\begin{equation*}\label{13}
\begin{aligned}
p = K_p\phi^{\frac{2\gamma}{\gamma-1}},
 \end{aligned}
\end{equation*}
\begin{equation*}\label{14}
c = \sqrt{-p_\tau} = K_c\phi^{\frac{\gamma+1}{\gamma-1}},
\end{equation*}
where $K_\tau$, $K_p$ and $K_c$ are positive constants given by
\begin{equation}\label{15}
\begin{aligned}
K_\tau := \left(\frac{2\sqrt{K\gamma}}{\gamma-1}\right)^{\frac{2}{\gamma-1}},\ \
K_p := K K_\tau^{-\gamma}, \ \ \text{and}\ \
K_c := \sqrt{K\gamma}K_\tau^{-\frac{\gamma+1}{2}}.
 \end{aligned}
\end{equation}
We also have
\begin{equation}\label{16}
\begin{aligned}
K_p=\frac{\gamma-1}{2\gamma}K_c \ \ \text{and}\ \ K_\tau K_c =\frac{\gamma-1}{2}.
 \end{aligned}
\end{equation}
In this paper, we always use $K$ with some subscripts to denote positive
constants. We will not notify the reader again if there is no ambiguity.

A direct calculation shows that the p-system have two characteristic speeds
$$\lambda_1=-\lambda_2=c.$$
The forward and backward characteristics are described by
\begin{equation}\label{17}
\begin{aligned}
\frac{dx}{dt}= c \ \ \text{and}\ \ \frac{dx}{dt}= -c,
 \end{aligned}
\end{equation}
and we denote the corresponding directional derivatives along them by
\begin{equation}\label{18}
\begin{aligned}
 ^\prime=\frac{\partial }{\partial t}+c\frac{\partial }{\partial x}
 \ \ \text{and}\ \
  ^\backprime=\frac{\partial }{\partial t}-c\frac{\partial }{\partial x}
 \end{aligned}
\end{equation}
respectively.  Furthermore, we denote the Riemann invariants are
\begin{equation}\label{19}
\begin{aligned}
{w}:=u+\phi \ \ \text{and}\ \  {z}:=u-\phi,
 \end{aligned}
\end{equation}
which satisfy
\begin{equation}\label{20.1}
\begin{aligned}
{w}^\prime:=-\frac{\a}{2(1+t)^\lambda}({z}+{w}),
 \end{aligned}
\end{equation}
and
\begin{equation}\label{20}
\begin{aligned}
{z}^\backprime:=-\frac{\a}{2(1+t)^\lambda}({z}+{w})
 \end{aligned}
\end{equation}
respectively.\\\\
\it Define
\begin{equation}\label{20.6}
\begin{aligned}
{A}={w}_x,\ \ {B}={z}_x,
 \end{aligned}
\end{equation}
then  we have:

\begin{lemma}\ \label{L2.2}
The smooth solutions  of (\ref{1})-(\ref{4}) satisfy
\begin{equation}\label{20.7}
\begin{aligned}
{A}^\prime= -\frac{\a}{2(1+t)^\lambda}({A}+{B})+K_d({A}{B}-{A}^2),
 \end{aligned}
\end{equation}
and
\begin{equation}\label{20.8}
\begin{aligned}
{B}^\backprime=-\frac{\a}{2(1+t)^\lambda}({A}+{B})+K_d({A}{B}-{B}^2),
 \end{aligned}
\end{equation}
where
$K_d = K_c\frac{\gamma+1}{2(\gamma-1)}\phi^{\frac{2}{\gamma-1}}.$
\end{lemma}

\textbf{Proof}. Take the partial derivative of both sides of equation (\ref{20.1}) with respect to $x$,
\begin{equation}\label{20.10}
\begin{aligned}
({w}^\prime)_x=-\frac{\a}{2(1+t)^\lambda}({z}+{w})_x,
 \end{aligned}
\end{equation}
then
\begin{equation}\label{20.11}
\begin{aligned}
({w}^\prime)_x=({w}_t+c{w}_x)_x={w}_{tx}+cw_{xx}+c_{x}{w}_{x}=({w}_x)^\prime+c_{x}{w}_{x}
=-\frac{\a}{2(1+t)^\lambda}({z}_x+{w}_x).
 \end{aligned}
\end{equation}
By (\ref{20.6}),
\begin{equation}\label{20.12}
\begin{aligned}
{A}^\prime+c_{x}{A}=-\frac{\a}{2(1+t)^\lambda}({A}+{B}).
 \end{aligned}
\end{equation}
Then  we obtain
\begin{equation}\label{20.13}
\begin{aligned}
c_x = \left(K_c\phi^{\frac{\gamma+1}{\gamma-1}}\right)_x
=K_c\frac{\gamma+1}{\gamma-1}\phi^{\frac{2}{\gamma-1}}\phi_x=2K_d\phi_x=K_d({A}-{B}),
 \end{aligned}
\end{equation}
where we used (\ref{19}), and denote $K_d=K_c\frac{\gamma+1}{2(\gamma-1)}\phi^{\frac{2}{\gamma-1}}$.
Thus (\ref{20.10}) changes to
\begin{equation}\label{20.15}
\begin{aligned}
{A}^\prime+K_d({A}-{B}){A}= -\frac{\a}{2(1+t)^\lambda}({A}+{B}),
 \end{aligned}
\end{equation}
and we get (\ref{20.7}). The calculation of (\ref{20.8}) is same.
This completes the proof. \\

\section{Singularity formation for compressible Euler equations when $\lambda\neq1$}
\subsection{¡°Decoupled¡± Riccati type equations }
\ \ \ \
To decouple $A$ and $B$ in (\ref{20.7}) and (\ref{20.8}) along the two directional derivatives,
we introduce two gradient variables
\begin{equation}\label{21}
\begin{aligned}
y:=\left(\phi^{\frac{\gamma+1}{2(\gamma-1)}}{A}-\frac{\a(\gamma-1)}{K_c(\gamma-3)(1+t)^\lambda}
\phi^{\frac{\gamma-3}{2(\gamma-1)}}\right)e^{\frac{\a(3\gamma-1)}{2(\gamma-3)(1-\lambda)}(1+t)^{1-\lambda}},
 \end{aligned}
\end{equation}
and
\begin{equation}\label{22}
\begin{aligned}
q:=\left(\phi^{\frac{\gamma+1}{2(\gamma-1)}}{B}-\frac{\a(\gamma-1)}{K_c(\gamma-3)(1+t)^\lambda}
\phi^{\frac{\gamma-3}{2(\gamma-1)}}\right)e^{\frac{\a(3\gamma-1)}{2(\gamma-3)(1-\lambda)}(1+t)^{1-\lambda}}.
 \end{aligned}
\end{equation}
We will show $y$ and $q$ satisfy the following Riccati equations:

\begin{lemma}\ \label{L3.1}
For $C^1$ solutions of (\ref{1})-(\ref{4}),  we have
\begin{equation}\label{23}
\begin{aligned}
y^\prime =a_0-a_2 y^2,
\end{aligned}
\end{equation}
\begin{equation}\label{24}
\begin{aligned}
q^\backprime =a_0-a_2 q^2,
\end{aligned}
\end{equation}
where
\begin{equation}\label{25}
\begin{aligned}
a_0=\frac{\lambda \a(\gamma-1)(\gamma-3)(1+t)^{\lambda-1}-\a^2(\gamma-1)^2}{K_c(\gamma-3)^2(1+t)^{2\lambda}}
\phi^{\frac{\gamma-3}{2(\gamma-1)}}e^{\frac{\a(3\gamma-1)}{2(\gamma-3)(1-\lambda)}(1+t)^{1-\lambda}},
\end{aligned}
\end{equation}
\begin{equation}\label{26}
\begin{aligned}
a_2=\frac{K_c(\gamma+1)}{2(\gamma-1)}
\phi^{-\frac{\gamma-3}{2(\gamma-1)}}e^{-\frac{\a(3\gamma-1)}{2(\gamma-3)(1-\lambda)}(1+t)^{1-\lambda}}>0.
\end{aligned}
\end{equation}
\end{lemma}

\textbf{Proof}. By (\ref{18}), (\ref{10}), (\ref{1}) and (\ref{19}),
\begin{equation}\label{26.1}
\begin{aligned}
\phi^\prime&=\phi_t+c\phi_x
=\left(\frac{2\sqrt{K\gamma}}{\gamma-1}\tau^{-\frac{\gamma-1}{2}}\right)_t+c\phi_x\\
&=-\sqrt{K\gamma}\tau^{-\frac{\gamma+1}{2}}\tau_t+c\phi_x\\
&=-cu_x+c\phi_x=-c({z}+\phi)_x+c\phi_x=-c{z}_x-c\phi_x+c\phi_x\\
&=-c{B}.
\end{aligned}
\end{equation}
Hence
\begin{equation}\label{26.2}
\begin{aligned}
{B}=-\frac{1}{c}\phi^\prime.
\end{aligned}
\end{equation}
Plugging (\ref{26.2}) into (\ref{20.7}),  we get
\begin{equation}\label{26.3}
\begin{aligned}
{A}^\prime= -\frac{\a}{2(1+t)^\lambda}({A}-\frac{1}{c}\phi^\prime)+K_d(-\frac{1}{c}\phi^\prime{A}-{A}^2).
 \end{aligned}
\end{equation}
 We move the terms including $\phi^\prime$ to the left hand side, then  we multiply by
$\phi^{\frac{\gamma+1}{2(\gamma-1)}}$ on both sides. After simplification,  we have
\begin{equation}\label{26.4}
\begin{aligned}
{A}^\prime\phi^{\frac{\gamma+1}{2(\gamma-1)}}
-\frac{\a}{2c(1+t)^\lambda}\phi^\prime\phi^{\frac{\gamma+1}{2(\gamma-1)}}
+\frac{K_d}{c}\phi^\prime{A}\phi^{\frac{\gamma+1}{2(\gamma-1)}}
= -\frac{\a}{2(1+t)^\lambda}{A}\phi^{\frac{\gamma+1}{2(\gamma-1)}}
-K_d{A}^2\phi^{\frac{\gamma+1}{2(\gamma-1)}}.
 \end{aligned}
\end{equation}
 The left hand side of (\ref{26.4}) is equal to
\begin{equation}\label{26.5}
\begin{aligned}
&{A}^\prime\phi^{\frac{\gamma+1}{2(\gamma-1)}}
-\frac{\a}{2c(1+t)^\lambda}\phi^\prime\phi^{\frac{\gamma+1}{2(\gamma-1)}}
+\frac{K_d}{c}\phi^\prime{A}\phi^{\frac{\gamma+1}{2(\gamma-1)}}\\
=&{A}^\prime\phi^{\frac{\gamma+1}{2(\gamma-1)}}
-\frac{\a}{2(1+t)^\lambda K_c\phi^{\frac{\gamma+1}{\gamma-1}}}\phi^\prime\phi^{\frac{\gamma+1}{2(\gamma-1)}}
+\frac{K_c\frac{\gamma+1}{2(\gamma-1)}\phi^{\frac{2}{\gamma-1}}}
{K_c\phi^{\frac{\gamma+1}{\gamma-1}}}\phi^\prime{A}\phi^{\frac{\gamma+1}{2(\gamma-1)}} \\
=&{A}^\prime\phi^{\frac{\gamma+1}{2(\gamma-1)}}
-\frac{\a}{2K_c(1+t)^\lambda}\phi^\prime\phi^{-\frac{\gamma+1}{2(\gamma-1)}}
+\frac{\gamma+1}{2(\gamma-1)}\phi^{\frac{-\gamma+3}{2(\gamma-1)}} \phi^\prime{A}
\\
=&\left({A}\phi^{\frac{\gamma+1}{2(\gamma-1)}}
-\frac{\a(\gamma-1)}{K_c(\gamma-3)}\phi^{\frac{\gamma-3}{2(\gamma-1)}}\frac{1}{(1+t)^\lambda}\right)^\prime
-\frac{\lambda \a(\gamma-1)}{K_c(\gamma-3)}\phi^{\frac{\gamma-3}{2(\gamma-1)}}(1+t)^{-\lambda-1}.
 \end{aligned}
\end{equation}
 We define a new variable $\tilde{y}$
\begin{equation}\label{26.6}
\begin{aligned}
\tilde{y}={A}\phi^{\frac{\gamma+1}{2(\gamma-1)}}
-\frac{\a(\gamma-1)}{K_c(\gamma-3)(1+t)^\lambda}\phi^{\frac{\gamma-3}{2(\gamma-1)}}.
 \end{aligned}
\end{equation}
So
\begin{equation}\label{26.7}
\begin{aligned}
{A}=\tilde{y}\phi^{-\frac{\gamma+1}{2(\gamma-1)}}
+\frac{\a(\gamma-1)}{K_c(\gamma-3)(1+t)^\lambda}\phi^{\frac{-4}{2(\gamma-1)}}.
 \end{aligned}
\end{equation}
 The  right hand side of (\ref{26.4}) is equal to
\begin{equation}\label{26.7.1}
\begin{aligned}
 &-\frac{\a}{2(1+t)^\lambda}{A}\phi^{\frac{\gamma+1}{2(\gamma-1)}}
-K_d{A}^2\phi^{\frac{\gamma+1}{2(\gamma-1)}}\\
=&-\frac{\a}{2(1+t)^\lambda}\phi^{\frac{\gamma+1}{2(\gamma-1)}}
\left(\tilde{y}\phi^{-\frac{\gamma+1}{2(\gamma-1)}}
+\frac{\a(\gamma-1)}{K_c(\gamma-3)(1+t)^\lambda}\phi^{\frac{-4}{2(\gamma-1)}}\right)\\
&-K_c\frac{\gamma+1}{2(\gamma-1)}\phi^{\frac{2}{\gamma-1}}\phi^{\frac{\gamma+1}{2(\gamma-1)}}
\left(\tilde{y}\phi^{-\frac{\gamma+1}{2(\gamma-1)}}
+\frac{\a(\gamma-1)}{K_c(\gamma-3)(1+t)^\lambda}\phi^{\frac{-4}{2(\gamma-1)}}\right)^2\\
=&-\frac{\a}{2(1+t)^\lambda}\tilde{y}
-\frac{\a^2(\gamma-1)}{2K_c(\gamma-3)(1+t)^{2\lambda}}\phi^{\frac{\gamma-3}{2(\gamma-1)}}
-\frac{\gamma+1}{2(\gamma-1)}K_c
\phi^{-\frac{\gamma-3}{2(\gamma-1)}}\tilde{y}^2\\
&-\frac{\a^2(\gamma^2-1)}{2K_c(\gamma-3)^2(1+t)^{2\lambda}}\phi^{\frac{\gamma-3}{2(\gamma-1)}}
-\frac{\a(\gamma+1)}{(\gamma-3)(1+t)^\lambda}\tilde{y}\\
=&-\frac{K_c(\gamma+1)}{2(\gamma-1)}\phi^{-\frac{\gamma-3}{2(\gamma-1)}}\tilde{y}^2
-\frac{\a(3\gamma-1)}{2(\gamma-3)(1+t)^\lambda}\tilde{y}
-\frac{\a^2(\gamma-1)^2}{K_c(\gamma-3)^2(1+t)^{2\lambda}}\phi^{\frac{\gamma-3}{2(\gamma-1)}}.
 \end{aligned}
\end{equation}
Then, (\ref{26.4}) changes to
\begin{equation}\label{26.8}
\begin{aligned}
\tilde{y}^\prime=& \frac{\lambda \a(\gamma-1)}{K_c(\gamma-3)}\phi^{\frac{\gamma-3}{2(\gamma-1)}}(1+t)^{-\lambda-1}
 -\frac{\a}{2(1+t)^\lambda}{A}\phi^{\frac{\gamma+1}{2(\gamma-1)}}
-K_d{A}^2\phi^{\frac{\gamma+1}{2(\gamma-1)}}\\
=&-\frac{K_c(\gamma+1)}{2(\gamma-1)}\phi^{-\frac{\gamma-3}{2(\gamma-1)}}\tilde{y}^2
-\frac{\a(3\gamma-1)}{2(\gamma-3)(1+t)^\lambda}\tilde{y}
-\frac{\a^2(\gamma-1)^2}{K_c(\gamma-3)^2(1+t)^{2\lambda}}\phi^{\frac{\gamma-3}{2(\gamma-1)}}\\
&+ \frac{\lambda \a(\gamma-1)}{K_c(\gamma-3)}\phi^{\frac{\gamma-3}{2(\gamma-1)}}(1+t)^{-\lambda-1}\\
=&-\frac{K_c(\gamma+1)}{2(\gamma-1)}\phi^{-\frac{\gamma-3}{2(\gamma-1)}}\tilde{y}^2
-\frac{\a(3\gamma-1)}{2(\gamma-3)(1+t)^\lambda}\tilde{y}
+\frac{\lambda \a(\gamma-1)(\gamma-3)(1+t)^{\lambda-1}-\a^2(\gamma-1)^2}
{K_c(\gamma-3)^2(1+t)^{2\lambda}}\phi^{\frac{\gamma-3}{2(\gamma-1)}}\\
:=&-\tilde{a}_2\tilde{y}^2-\tilde{a}_1\tilde{y}+\tilde{a}_0,
 \end{aligned}
\end{equation}
where
\begin{equation}\label{26.8.1}
\begin{aligned}
&\tilde{a}_2=\frac{K_c(\gamma+1)}{2(\gamma-1)}\phi^{-\frac{\gamma-3}{2(\gamma-1)}},\\
&\tilde{a}_1=\frac{\a(3\gamma-1)}{2(\gamma-3)(1+t)^\lambda},\\
&\tilde{a}_0=\frac{\lambda \a(\gamma-1)(\gamma-3)(1+t)^{\lambda-1}-\a^2(\gamma-1)^2}
{K_c(\gamma-3)^2(1+t)^{2\lambda}}\phi^{\frac{\gamma-3}{2(\gamma-1)}}.
 \end{aligned}
\end{equation}
Then  we do one more simplification by multiplying
\begin{equation}\label{26.9}
\begin{aligned}
\tilde{\mu}=e^{\frac{\a(3\gamma-1)}{2(\gamma-3)(1-\lambda)}(1+t)^{1-\lambda}}
 \end{aligned}
\end{equation}
on (\ref{26.8}). In fact, it is easy to check that
\begin{equation}\label{26.10}
\begin{aligned}
\tilde{\mu}^\prime=\tilde{\mu}_t+c\tilde{\mu}_x=\tilde{a}_1\tilde{\mu}.
 \end{aligned}
\end{equation}
Then  we denote
\begin{equation}\label{26.11}
\begin{aligned}
y=\tilde{\mu}\tilde{y}.
 \end{aligned}
\end{equation}
Hence (\ref{26.8}) changes to
\begin{equation}\label{26.12}
\begin{aligned}
y^\prime=a_0-a_2y^2.
 \end{aligned}
\end{equation}
where
\begin{equation}\label{26.13}
\begin{aligned}
a_0=\tilde{\mu}\tilde{a}_0, \ \ a_2=\frac{\tilde{a}_2}{\tilde{\mu}}.
 \end{aligned}
\end{equation}

Similarly, we prove (\ref{24}).
Firstly, we have
\begin{equation}\label{2601}
\begin{aligned}
\phi^\backprime&=\phi_t-c\phi_x
=\left(\frac{2\sqrt{K\gamma}}{\gamma-1}\tau^{-\frac{\gamma-1}{2}}\right)_t-c\phi_x\\
&=-\sqrt{K\gamma}\tau^{-\frac{\gamma+1}{2}}\tau_t-c\phi_x\\
&=-cu_x-c\phi_x=-c({w}-\phi)_x-c\phi_x=-c{w}_x+c\phi_x-c\phi_x\\
&=-c{A},
\end{aligned}
\end{equation}
and
\begin{equation}\label{2603}
\begin{aligned}
{B}^\backprime= -\frac{\a}{2(1+t)^\lambda}({B}-\frac{1}{c}\phi^\backprime)+K_d(-\frac{1}{c}\phi^\backprime{B}-{B}^2).
 \end{aligned}
\end{equation}
By calculating, we define a new variable
\begin{equation}\label{2606}
\begin{aligned}
\tilde{q}={B}\phi^{\frac{\gamma+1}{2(\gamma-1)}}
-\frac{\a(\gamma-1)}{K_c(\gamma-3)(1+t)^\lambda}\phi^{\frac{\gamma-3}{2(\gamma-1)}},
 \end{aligned}
\end{equation}
and then
\begin{equation}\label{2607}
\begin{aligned}
{B}=\tilde{q}\phi^{-\frac{\gamma+1}{2(\gamma-1)}}
+\frac{\a(\gamma-1)}{K_c(\gamma-3)(1+t)^\lambda}\phi^{\frac{-4}{2(\gamma-1)}}.
 \end{aligned}
\end{equation}
We can get that $\tilde{q}$ satisfies
\begin{equation}\label{2608}
\begin{aligned}
\tilde{q}^\backprime=& \frac{\lambda \a(\gamma-1)}{K_c(\gamma-3)}\phi^{\frac{\gamma-3}{2(\gamma-1)}}(1+t)^{-\lambda-1}
 -\frac{\a}{2(1+t)^\lambda}{B}\phi^{\frac{\gamma+1}{2(\gamma-1)}}
-K_d{B}^2\phi^{\frac{\gamma+1}{2(\gamma-1)}}\\
=&-\frac{K_c(\gamma+1)}{2(\gamma-1)}\phi^{-\frac{\gamma-3}{2(\gamma-1)}}\tilde{q}^2
-\frac{\a(3\gamma-1)}{2(\gamma-3)(1+t)^\lambda}\tilde{q}
-\frac{\a^2(\gamma-1)^2}{K_c(\gamma-3)^2(1+t)^{2\lambda}}\phi^{\frac{\gamma-3}{2(\gamma-1)}}\\
&+ \frac{\lambda \a(\gamma-1)}{K_c(\gamma-3)}\phi^{\frac{\gamma-3}{2(\gamma-1)}}(1+t)^{-\lambda-1}\\
=&-\frac{K_c(\gamma+1)}{2(\gamma-1)}\phi^{-\frac{\gamma-3}{2(\gamma-1)}}\tilde{q}^2
-\frac{\a(3\gamma-1)}{2(\gamma-3)(1+t)^\lambda}\tilde{q}
+\frac{\lambda \a(\gamma-1)(\gamma-3)(1+t)^{\lambda-1}-\a^2(\gamma-1)^2}
{K_c(\gamma-3)^2(1+t)^{2\lambda}}\phi^{\frac{\gamma-3}{2(\gamma-1)}}\\
:=&-\tilde{a}_2\tilde{q}^2-\tilde{a}_1\tilde{q}+\tilde{a}_0,
 \end{aligned}
\end{equation}
where $\tilde{a}_2,$ $\tilde{a}_1,$ $\tilde{a}_0$ are defined in  (\ref{26.8.1}).
Again, if we denote $q=\tilde{\mu}\tilde{q}$,
then (\ref{2608}) changes to
\begin{equation}\label{26.012}
\begin{aligned}
q^\backprime=a_0-a_2q^2.
 \end{aligned}
\end{equation}
This completes the proof.\\

\subsection{Uniform lower bounds on density for $1<\gamma<3$}
\ \ \ \
 In this subsection, the estimates of density's lower bounds will be considered for adiabatic exponent $1<\gamma<3$.
 Firstly, noticing the definition of $a_0$ in (\ref{25}), we get
\begin{equation}\label{B1}
1<\gamma<3,~~ \ld\g {{\a(\gamma-1)}\over{\gamma-3}}
\end{equation}
make $a_0\l 0$. From Lemma \ref{L3.1}, we have
\begin{lemma}\ \label{L3.2}
 Suppose (\ref{B1}) satisfies. Then for any $C^1$ solutions of  (\ref{1})-(\ref{4}), we have a priori bounds
\begin{equation}\label{27}
\begin{aligned}
y(x, t)\leq \max\{1, \sup_x(y(x,0))\} =: Y,
\end{aligned}
\end{equation}
and
\begin{equation}\label{28}
\begin{aligned}
q(x,t)\leq  \max \{1, \sup_x(q(x,0))\} =: Q.
\end{aligned}
\end{equation}
\end{lemma}
The following lemma provide us a time-dependent lower bound on density for arbitrary
classical solutions.

\begin{lemma}\ \label{L3.3}
 Assume (\ref{B1}) fulfil and (\ref{1})-(\ref{4}) has a $C^1$ solution $(\tau, u)(x,t)$. Then there exists a positive constant $T>0$, such that
\begin{equation}\label{29}
\begin{aligned}
\rho(x,t)=\tau^{-1}(x,t)\geq K_0t^{-\frac{4}{3-\gamma}}e^{\frac{-2\a(3\gamma-1)(1+t)^{1-\lambda}}{(3-\gamma)^2(1-\lambda)}},\ \ \text{for any} \ \ t>T,
\end{aligned}
\end{equation}
where $K_0=\left[\left(\frac{2\sqrt{K\gamma}}{\gamma-1}\right)^{-\frac{\gamma-3}{2(\gamma-1)}}(Y+Q)
\frac{3-\gamma}{2(\gamma-1)}K_c\right]^{-\frac{4}{3-\gamma}}.$
\end{lemma}

\textbf{Proof}. From (\ref{20.1}), (\ref{21}) and Lemma \ref{L3.2},  we have
\begin{equation}\label{30}
\begin{aligned}
{w}_t&=-c{w}_x-\frac{\a}{2(1+t)^\lambda}({z}+{w})\\
&=-K_c\phi^{\frac{\gamma+1}{\gamma-1}}
\left( ye^{-\frac{\a(3\gamma-1)(1+t)^{1-\lambda}}{2(\gamma-3)(1-\lambda)}}\phi^{-\frac{\gamma+1}{2(\gamma-1)}}
+\frac{\a(\gamma-1)}{K_c(\gamma-3)(1+t)^\lambda}
\phi^{\frac{-2}{\gamma-1}} \right)-\frac{\a}{2(1+t)^\lambda}({w}+{z})\\
&=-K_c\phi^{\frac{\gamma+1}{2(\gamma-1)}}
e^{-\frac{\a(3\gamma-1)(1+t)^{1-\lambda}}{2(\gamma-3)(1-\lambda)}}y
-\frac{\a(\gamma-1)}{(\gamma-3)(1+t)^\lambda}\phi-\frac{\a}{2(1+t)^\lambda}({w}+{z})\\
&\geq-K_c\phi^{\frac{\gamma+1}{2(\gamma-1)}}e^{-\frac{\a(3\gamma-1)(1+t)^{1-\lambda}}{2(\gamma-3)(1-\lambda)}}Y
-\frac{\a(\gamma-1)}{(\gamma-3)(1+t)^\lambda}\phi-\frac{\a}{2(1+t)^\lambda}({w}+{z}).
\end{aligned}
\end{equation}
 Similarly, from (\ref{20}), (\ref{22})  and Lemma \ref{L3.2},  we have
\begin{equation}\label{31}
\begin{aligned}
{z}_t&=c{z}_x-\frac{\a}{2(1+t)^\lambda}({z}+{w})\\
&=K_c\phi^{\frac{\gamma+1}{\gamma-1}}
\left(qe^{-\frac{\a(3\gamma-1)(1+t)^{1-\lambda}}{2(\gamma-3)(1-\lambda)}}\phi^{-\frac{\gamma+1}{2(\gamma-1)}}
+\frac{\a(\gamma-1)}{K_c(\gamma-3)(1+t)^\lambda}
\phi^{\frac{-2}{\gamma-1}} \right)-\frac{\a}{2(1+t)^\lambda}({w}+{z})\\
&=K_c\phi^{\frac{\gamma+1}{2(\gamma-1)}}e^{-\frac{\a(3\gamma-1)(1+t)^{1-\lambda}}{2(\gamma-3)(1-\lambda)}}q
+\frac{\a(\gamma-1)}{(\gamma-3)(1+t)^\lambda}\phi-\frac{\a}{2(1+t)^\lambda}({w}+{z})\\
&\leq K_c\phi^{\frac{\gamma+1}{2(\gamma-1)}}e^{-\frac{\a(3\gamma-1)(1+t)^{1-\lambda}}{2(\gamma-3)(1-\lambda)}}Q
+\frac{\a(\gamma-1)}{(\gamma-3)(1+t)^\lambda}\phi-\frac{\a}{2(1+t)^\lambda}({w}+{z}).
\end{aligned}
\end{equation}
Therefore,
\begin{equation}\label{32}
\begin{aligned}
\phi_t&=\frac{1}{2}({w}_t-{z}_t)\\
&\geq-\frac{1}{2}K_c\phi^{\frac{\gamma+1}{2(\gamma-1)}}e^{-\frac{\a(3\gamma-1)(1+t)^{1-\lambda}}{2(\gamma-3)(1-\lambda)}}(Y+Q)
-\frac{\a(\gamma-1)}{(\gamma-3)(1+t)^\lambda}\phi\\
&\geq-\frac{1}{2}K_c\phi^{\frac{\gamma+1}{2(\gamma-1)}}e^{\frac{\a(3\gamma-1)(1+t)^{1-\lambda}}{2(3-\gamma)(1-\lambda)}}(Y+Q),
\end{aligned}
\end{equation}
where we have used $1<\gamma<3$.
Dividing the above inequality by $\phi^{\frac{\gamma+1}{2(\gamma-1)}}$, then integrating both sides on $t$,   we obtain
\begin{equation}\label{33}
\begin{aligned}
&\frac{2(\gamma-1)}{\gamma-3}(\phi(x,t))^{\frac{\gamma-3}{2(\gamma-1)}}-\frac{2(\gamma-1)}{\gamma-3}(\phi(x,0))^{\frac{\gamma-3}{2(\gamma-1)}}\\
\geq&-\frac{1}{2}K_c(Y+Q)\int_0^te^{\frac{\a(3\gamma-1)(1+s)^{1-\lambda}}{2(3-\gamma)(1-\lambda)}}ds\\
\geq&-\frac{1}{2}K_c(Y+Q)\int_0^te^{\frac{\a(3\gamma-1)(1+t)^{1-\lambda}}{2(3-\gamma)(1-\lambda)}}ds\\
=&-\frac{1}{2}K_c(Y+Q)te^{\frac{\a(3\gamma-1)(1+t)^{1-\lambda}}{2(3-\gamma)(1-\lambda)}}.
\end{aligned}
\end{equation}
In above calculation, we have used the function $e^{{(1+t)^{1-\ld}}\over {1-\ld}}$ is increasing with $t$ for any parameter $\ld\in R$.
Then
\begin{equation}\label{33.1}
\begin{aligned}
(\phi(x,t))^{\frac{\gamma-3}{2(\gamma-1)}}&\leq(\phi(x,0))^{\frac{\gamma-3}{2(\gamma-1)}}
+\frac{3-\gamma}{4(\gamma-1)}K_cte^{\frac{\a(3\gamma-1)(1+t)^{1-\lambda}}{2(3-\gamma)(1-\lambda)}}(Y+Q).
\end{aligned}
\end{equation}
There exists a positive constant $T$, when $t>T$,
\begin{equation}\label{33.2}
\begin{aligned}
(\phi(x,t))^{\frac{\gamma-3}{2(\gamma-1)}}
\leq\frac{3-\gamma}{2(\gamma-1)}K_cte^{\frac{\a(3\gamma-1)(1+t)^{1-\lambda}}{2(3-\gamma)(1-\lambda)}}(Y+Q),\ \ t>T.
\end{aligned}
\end{equation}
Therefore,
\begin{equation}\label{33.3}
\begin{aligned}
\rho(x,t)&\geq\left[\left(\frac{2\sqrt{K\gamma}}{\gamma-1}\right)^{-\frac{\gamma-3}{2(\gamma-1)}}
(Y+Q)\frac{3-\gamma}{2(\gamma-1)}K_c\right]^{-\frac{4}{3-\gamma}}
t^{-\frac{4}{3-\gamma}}e^{\frac{-2\a(3\gamma-1)(1+t)^{1-\lambda}}{(3-\gamma)^2(1-\lambda)}}\\
&:=K_0t^{-\frac{4}{3-\gamma}}e^{\frac{-2\a(3\gamma-1)(1+t)^{1-\lambda}}{(3-\gamma)^2(1-\lambda)}},\ \ t>T.
\end{aligned}
\end{equation}
This completes the proof.\\

\subsection{Singularity formation mechanism}
\ \ \ \
 In this subsection, we will investigate the decoupled ODEs (\ref{23}) and (\ref{24})
 to gain some sufficient conditions for singularity formation of the system (\ref{1})-(\ref{4}) in finite time.
 We first consider the blow up mechanism of ODEs (\ref{23}) and (\ref{24}) alone the forward and backward characteristic.
 And we give the analysis of the following two cases.\\

  \it \textbf{1)} \textit{$a_0\l 0$ and $y(x_0,0)=y(x(0),0)\l 0$ for some $x_0\in R$.} \\

\it In this case, (\ref{23}) changes into
\begin{equation}\label{36}
\begin{aligned}
y^\prime =a_0-a_2 y^2\l-a_2 y^2\leq0.
\end{aligned}
\end{equation}
Integrating both sides of the above inequality with respect to $t$ along the forward characteristic,
we have
\begin{equation}\label{38}
\begin{aligned}
0>y^{-1}(x(t),t)\geq y^{-1}(x_0,0)+\int_0^ta_2(x(s),s)ds,
\end{aligned}
\end{equation}
and
\begin{equation}\label{39}
\begin{aligned}
y(x(t),t)\leq\frac{1}{y^{-1}(x_0,0)+\int_0^ta_2(x(s),s)ds}.
\end{aligned}
\end{equation}
From (\ref{39}),  we see that, for some positive constant $\beta>0$,
if $-\beta<y^{-1}(x_0,0)<0$ and $\int_0^\infty a_2(x(s),s)ds\geq \beta$,
then there must exist a finite time $t_*$ such that $y(x(t),t)$ blow up before $t^*$.\\\\

\it \textbf{2)} \textit{$a_0>0$ and $y(x_0,0)<-\disp\max_{(x,t)}\sqrt{a_0\over a_2}:=-N$ for some $x_0\in R$.}\\

\it In this case, (\ref{23}) turns into
\begin{equation}\label{36}
\begin{aligned}
y^\prime =-a_2\left(y-\sqrt{{a_0}\over{a_2}}\right)\left(y+\sqrt{{a_0}\over{a_2}}\right),
\end{aligned}
\end{equation}
and there exists a $\varepsilon> 0$ such that
$$y(x_0, 0) < -(1 + \varepsilon) N .$$
Along the forward characteristic,
recall (\ref{23}), we have
\begin{equation}\label{36}
y'< 0 ~~{\rm{and}}~~ y(x(t),t) < -(1 + \varepsilon) N, ~~{\rm for~ any}~ t\geq 0,
\end{equation}
which together with the definition of $N$ implies
$$
a_0-a_2 {y^2(x(t),t)
\over(1 + \varepsilon)^2} < 0, ~~{\rm{for~ all}}~ t\geq 0.
$$
Hence by (\ref{23}) and $a_2 > 0$,
\begin{equation}\label{}
\begin{aligned}
y^\prime <\left(-1+{1\over {(1+\varepsilon)^2}}\right)a_2y^2.
\end{aligned}
\end{equation}
Like the calculations in $(\ref{39})$, we have
\begin{equation}\label{3.2.7}
\begin{aligned}
y(x(t),t)<\left({y^{-1}(x_0,0)+\left(1-{1\over{(1+\varepsilon)^2}}\right)\disp\int_0^ta_2(x(s),s)ds}\right)^{-1},
\end{aligned}
\end{equation}
where the integral is along the forward characteristic.
Again, we see that
if $y(x_0,0)<-N<0$ and $\int_0^\infty a_2(x(s),s)ds=+\i$,
then there must exist a finite time $t_*$ such that $y(x(t),t)$ blow up before $t^*$.\\
\it To conclude above analysis, we give the following Lemma:\\
\begin{lemma}\label{L3.4}
Consider the differential equation
$$
y'=y_t+cy_x=a_0-a_2y^2,
$$
where $c, a_0, a_2$ are some functions which may depend on $x$ and $t$. We suppose $a_2>0$ for any $(x,t)\in R\times R^+$.
Denote $x(t)$ be the characteristic line start from $(x_0,0)$, i.e. $x(t)$ satisfies the ODE
$$
{{dx(t)}\over {dt}}=c(x(t),t),~~~x(0)=x_0.
$$
If one of the following two conditions\\
\it \textbf{1)}  $a_0\l 0$ and there exists one point $x_0\in R$ such that $\disp y(x_0,0)< -\left(\int_0^\i a_2(x(s),s)ds\right)^{-1}$;\\
\it \textbf{2)}  $a_0>0, ~ \disp\int_0^\i a_2(x(s),s)ds=+\i$   and there exists one point $x_0\in R$ such that
\begin{equation}
y(x_0,0)\l -\max_{(x,t)}\sqrt{{a_0}\over{a_2}}
\end{equation}
satisfies. Then there must exist a finite time $t_*$ such that
$$
\lim_{t\to t_*-}y(x(t),t)=-\i.
$$
\end{lemma}

\it From Lemma \ref{L3.4}, in order to consider the blow up mechanism, it is important to gain the integration $\int_0^\i a_2(x(s),s)ds$.
 Base on the expression of $a_2$, the adiabatic exponent need to be separated into two cases: $\gamma>3$ and $1<\gamma<3$,
 which rely on the upper and lower bounds of density respectively.\\\\

 \it \textit{\textbf{Case I:} For $\gamma>3$.}\\
\it From the definition of $a_0$ in (\ref{25}), We claim that $a_0\l 0$ is equivalent to
\begin{equation}\label{3.2.1}
\ld (\gamma-3)\l \a(\gamma-1)(1+t)^{1-\ld}
\end{equation}
after a simple calculation. Therefore,
\begin{equation}\label{3.2.2}
\ld< \min\big\{1~~ ,{\a(\gamma-1)\over {\gamma-3}}\big\}
\end{equation}
makes $a_0\l 0$.\\
\it Using (\ref{26}) and Theorem \ref{t2.1}, we calculate
\begin{equation}\label{41}
\begin{aligned}
\int_0^\infty a_2(x(s),s)ds
&=\int_0^\infty\frac{K_c(\gamma+1)}{2(\gamma-1)}
\phi^{-\frac{\gamma-3}{2(\gamma-1)}}e^{-\frac{\a(3\gamma-1)}{2(\gamma-3)(1-\lambda)}(1+s)^{1-\lambda}}ds\\
&=\int_0^\infty\frac{K_c(\gamma+1)}{2(\gamma-1)}\left(\frac{2\sqrt{K\gamma}}{\gamma-1}\rho^{\frac{\gamma-1}{2}}\right)^
{-\frac{\gamma-3}{2(\gamma-1)}}e^{-\frac{\a(3\gamma-1)}{2(\gamma-3)(1-\lambda)}(1+s)^{1-\lambda}}ds\\
&=\frac{K_c(\gamma+1)}{2(\gamma-1)}\left(\frac{2\sqrt{K\gamma}}{\gamma-1}\right)^{-\frac{\gamma-3}{2(\gamma-1)}}
\int_0^\infty\rho^{\frac{3-\gamma}{4}}e^{-\frac{\a(3\gamma-1)}{2(\gamma-3)(1-\lambda)}(1+s)^{1-\lambda}}ds\\
&\geq\frac{K_c(\gamma+1)}{2(\gamma-1)}\left(\frac{2\sqrt{K\gamma}}{\gamma-1}\right)^{-\frac{\gamma-3}{2(\gamma-1)}}
\tilde{C}_0^{\frac{3-\gamma}{4}}\int_0^\infty e^{-\frac{\a(3\gamma-1)}{2(\gamma-3)(1-\lambda)}(1+s)^{1-\lambda}}ds,\\
&:= K_1I_1.
\end{aligned}
\end{equation}
where $K_1=\disp\frac{K_c(\gamma+1)}{2(\gamma-1)}\left(\frac{2\sqrt{K\gamma}}{\gamma-1}\right)^{-\frac{\gamma-3}{2(\gamma-1)}}
\tilde{C}_0^{\frac{3-\gamma}{4}}, ~I_1=\disp\int_0^\infty e^{-\frac{\a(3\gamma-1)}{2(\gamma-3)(1-\lambda)}(1+s)^{1-\lambda}}ds.$\\
\it Now, we analysis the lower bounds of $I_1$. When $0\l\ld< 1~~ \rm{and}~~ \ld\l {\a(\gamma-1)\over {\gamma-3}}$,
\begin{equation}\label{3.2.3}
\begin{aligned}
I_1&>\int_{0}^\i (1+s)^{-\ld}e^{-\frac{\a(3\gamma-1)}{2(\gamma-3)(1-\lambda)}(1+s)^{1-\lambda}}ds\\
&={{2(\gamma-3)}\over{\a(3\gamma-1)}}e^{-\frac{\a(3\gamma-1)}{2(\gamma-3)(1-\lambda)}}.
\end{aligned}
\end{equation}
While when $\ld\l 0$, $I_1$ is integrable. Let
\begin{equation}K_2=\left\{ \begin{aligned}{{2(\gamma-3)}\over{\a(3\gamma-1)}}e^{-\frac{\a(3\gamma-1)}{2(\gamma-3)(1-\lambda)}}, ~~  \text{for} \ 0\l\ld< 1,\\
\int_0^{\infty}e^{-\frac{\a(3\gamma-1)}{2(\gamma-3)(1-\lambda)}(1+s)^{1-\lambda}}ds, ~~  \text{for} \ \ld\l 0. \end{aligned}\right.\end{equation}Then
\begin{equation}\label{3.2.5}
\begin{aligned}
\int_0^\infty a_2(x(s),s)ds>K_1K_2>0
\end{aligned}
\end{equation}
according to (\ref{41}). If we suppose the initial data $y(x_0, 0)<-{1\over{K_1K_2}}$,
 then $y(x(t), t)$ must blow up in finite time.\\
\it On the other hand, $a_0> 0$ is equivalent to
\begin{equation}\label{3.2.8}
\ld (\gamma-3)> \a(\gamma-1)(1+t)^{1-\ld}.
\end{equation}
Therefore
\begin{equation}\label{B2}
\ld>\max\big\{1,{{\a(\gamma-1)}\over {\gamma-3}}\big\}
\end{equation}
makes $a_0> 0$ and
\begin{equation}\label{3.2.9}
I_1=\int_0^\i e^{-{{\a(3\gamma-1)}\over{2(\gamma-3)(1-\ld)}}(1+s)^{1-\ld}}ds>\int_0^\i ds=+\i.
\end{equation}
While
\begin{equation}\label{3.2.10}
\begin{aligned}
{a_0\over a_2}&=\frac{2\a(\gamma-1)^2[\lambda(\gamma-3) (1+t)^{\lambda-1}-\a(\gamma-1)]}{K_c^2(\gamma-3)^2(1+t)^{2\lambda}(\gamma+1)}
\phi^{\frac{\gamma-3}{\gamma-1}}e^{\frac{\a(3\gamma-1)}{(\gamma-3)(1-\lambda)}(1+t)^{1-\lambda}},\\
&\l K_3\tilde{C}_0^{\frac{\gamma-3}{2}}\lambda(\gamma-3) (1+t)^{-\lambda-1}e^{\frac{\a(3\gamma-1)}{(\gamma-3)(1-\lambda)}(1+t)^{1-\lambda}}\\
&\l K_3\tilde{C}_0^{\frac{\gamma-3}{2}}\lambda(\gamma-3):=K_4^2,
\end{aligned}
\end{equation}
where $K_3=\frac{2\a(\gamma-1)^2}{K_c^2(\gamma-3)^2(\gamma+1)}\big({2\sqrt{K\gamma}\over{\gamma-1}}\big)^{{\gamma-3}\over{\gamma-1}}>0,~ K_4>0$.
If we suppose the initial data $y(x_0, 0)<-K_4$, then $y(x(t), t)$ must blow up in finite time.\\\\
\it To sum up, suppose $C_0$ is the uniform upper bound of initial data $|u_0(x)|$ and $\tau_0^{-1}(x)$, we denote the constants
\begin{equation}\label{3.2.13}
\begin{aligned}
&\tilde{C}_0=\max\{C_0+C_0^\theta, (C_0+C_0^\theta)^{1\over \theta}\},~~~\theta={{\gamma-1}\over 2},\\
& K_1=\disp\frac{K_c(\gamma+1)}{2(\gamma-1)}\left(\frac{2\sqrt{K\gamma}}{\gamma-1}\right)
^{-\frac{\gamma-3}{2(\gamma-1)}}\tilde{C}_0^{\frac{3-\gamma}{4}},\\
&K_2=\left\{ \begin{aligned}{{2(\gamma-3)}\over{\a(3\gamma-1)}}e^{-\frac{\a(3\gamma-1)}{2(\gamma-3)(1-\lambda)}}, ~~  \text{for} \ 0\l\ld< 1,\\
\int_0^{\infty}e^{-\frac{\a(3\gamma-1)}{2(\gamma-3)(1-\lambda)}(1+s)^{1-\lambda}}ds, ~~  \text{for} \ \ld\l 0, \end{aligned}\right.\\
&K_3=\frac{2\a(\gamma-1)^2}{K_c^2(\gamma-3)^2(\gamma+1)}\big({2\sqrt{K\gamma}\over{\gamma-1}}\big)^{{\gamma-3}\over{\gamma-1}}>0,
~~K_4=\sqrt{K_3\tilde{C}_0^{\frac{\gamma-3}{2}}\lambda(\gamma-3)},\\
&N=\left\{ \begin{aligned}{1\over {K_1K_2}},~~\text{for} ~~\ld<\min\big\{1, {{\alpha(\gamma-1)}\over{\gamma-3}} \big\}, \\
K_4,~~\text{for} ~~\ld>\max\big\{1, {{\alpha(\gamma-1)}\over{\gamma-3}}\big\}.\end{aligned} \right.
\end{aligned}
\end{equation}
 Next, we give the singularity formation theorem for $\gamma>3$ in the following:
\begin{theorem}\label{t3.1}(Singularity formation for $\gamma>3$)\\
\it Suppose the initial data $(\tau_0, u_0)(x)\in C^1(R)$ and there exists a positive constant $C_0$ such that
$$|u_0(x)|\l C_0,~ 0< \tau_0^{-1}\l C_0. $$
Let $\ld< \min\{1, {{\alpha(\gamma-1)}\over{\gamma-3}} \}$ (or $\ld>\max\{1, {{\alpha(\gamma-1)}\over{\gamma-3}}\}$ ).
 Assume there exists one point $x_0$ such that $y(x_0,0)\l -N$ or $q(x_0,0)\l -N$ , i.e.
\begin{equation}\label{3.2.11}
\begin{aligned}
&u_x(x_0,0)+\phi_x(x_0,0)<\tilde{K}_1(\phi(x_0,0))^{\frac{-2}{\gamma-1}}-\tilde{K}_2(\phi(x_0,0))^{-\frac{\gamma+1}{2(\gamma-1)}},
\end{aligned}
\end{equation}
or
\begin{equation}\label{3.2.12}
\begin{aligned}
&u_x(x_0,0)-\phi_x(x_0,0)<\tilde{K}_1(\phi(x_0,0))^{\frac{-2}{\gamma-1}}-\tilde{K}_2(\phi(x_0,0))^{-\frac{\gamma+1}{2(\gamma-1)}},
\end{aligned}
\end{equation}
where
$\tilde{K}_1=\frac{\a(\gamma-1)}{K_c(\gamma-3) },
~\tilde{K}_2=Ne^{-\frac{\a(3\gamma-1)}{2(\gamma-3)(1-\lambda)} },~\phi =\frac{2\sqrt{K\gamma}}{\gamma-1}\tau^{-\frac{\gamma-1}{2}}$. Then $u_x$ and/or $\tau_x$ must blow up in finite time.
\end{theorem}

\it\textit{\textbf{Case II:} For $1<\gamma<3$.}\\
 \it From Lemma \ref{L3.3},  we have
\begin{equation}\label{42}
\begin{aligned}
(\rho(x,t))^{\frac{3-\gamma}{4}}\geq \left(K_0t^{-\frac{4}{3-\gamma}}e^{\frac{-2\a(3\gamma-1)(1+t)^{1-\lambda}}{(3-\gamma)^2(1-\lambda)}}\right)^{\frac{3-\gamma}{4}},\ \ t>T,
\end{aligned}
\end{equation}
then  there exists a positive constant $T$ such that
\begin{equation}\label{43}
\begin{aligned}
&\ \ \ \ \int_0^\infty a_2(x(s),s)ds\\
&=\frac{K_c(\gamma+1)}{2(\gamma-1)}\left(\frac{2\sqrt{K\gamma}}{\gamma-1}\right)^{-\frac{\gamma-3}{2(\gamma-1)}}
\left(\int_0^T\rho^{\frac{3-\gamma}{4}}e^{-\frac{\a(3\gamma-1)}{2(\gamma-3)(1-\lambda)}(1+s)^{1-\lambda}}ds
+\int_T^\infty\rho^{\frac{3-\gamma}{4}}e^{-\frac{\a(3\gamma-1)}{2(\gamma-3)(1-\lambda)}(1+s)^{1-\lambda}}\right)\\
&\geq\frac{K_c(\gamma+1)}{2(\gamma-1)}\left(\frac{2\sqrt{K\gamma}}{\gamma-1}\right)^{-\frac{\gamma-3}{2(\gamma-1)}}
\int_T^\infty\rho^{\frac{3-\gamma}{4}}e^{-\frac{\a(3\gamma-1)}{2(\gamma-3)(1-\lambda)}(1+s)^{1-\lambda}}ds\\
&\geq\frac{K_c(\gamma+1)}{2(\gamma-1)}\left(\frac{2\sqrt{K\gamma}}{\gamma-1}\right)^{-\frac{\gamma-3}{2(\gamma-1)}}
\int_T^\infty\left(K_0s^{-\frac{4}{3-\gamma}}e^{\frac{-2\a(3\gamma-1)(1+s)^{1-\lambda}}{(3-\gamma)^2(1-\lambda)}}\right)
^{\frac{3-\gamma}{4}}e^{-\frac{\a(3\gamma-1)}{2(\gamma-3)(1-\lambda)}(1+s)^{1-\lambda}}ds\\
&=\frac{K_c(\gamma+1)}{2(\gamma-1)}\left(\frac{2\sqrt{K\gamma}}{\gamma-1}\right)^{-\frac{\gamma-3}{2(\gamma-1)}}K_0^{\frac{3-\gamma}{4}}
\int_T^\infty s^{-1}ds=\infty.\ \ \ \
\end{aligned}
\end{equation}
Using Case 1) in Lemma 3.4, when $\ld\g \frac{\a(\gamma-1)}{(\gamma-3)}$, we have:
\begin{theorem} (Singularity formation for $1<\gamma<3$)\\
\it Let the initial data $(\tau_0, u_0)(x)\in C^1(R)$ and there exists a positive constant $C_0$ such that
$$|u_0(x)|\l C_0,~ 0< \tau_0^{-1}\l C_0. $$
 Suppose $\ld \g {{\a(\gamma-1)}\over{\gamma-3}}$ and there exists one point $x_0$ such that $y(x_0,0)\l 0$ or $q(x_0,0)\l 0$ , i.e.
\begin{equation}\label{35}
\begin{aligned}
u_x(x_0,0)+\phi_x(x_0,0)<-\frac{\a(\gamma-1)}{K_c(3-\gamma)}(\phi(x_0,0))^{\frac{-2}{\gamma-1}},
\end{aligned}
\end{equation}
or
\begin{equation}\label{35.0}
\begin{aligned}
u_x(x_0,0)-\phi_x(x_0,0)<-\frac{\a(\gamma-1)}{K_c(3-\gamma)}(\phi(x_0,0))^{\frac{-2}{\gamma-1}},
\end{aligned}
\end{equation}
then $u_x$ and/or $\phi_x$ blow up in finite time.
\end{theorem}

\section{Singularity formation for compressible Euler equations when $\lambda=1$}

\subsection{Decoupled ordinary differential equation }
\ \ \ \
Like the calculations in Section 3, to decouple $A$ and $B$ in (\ref{20.7}) and (\ref{20.8}) along the two directional derivatives when $\lambda=1$, we introduce two gradient variables
\begin{equation}\label{4.1}
\begin{aligned}
{y_1}:=\left(\phi^{\frac{\gamma+1}{2(\gamma-1)}}{A}-\frac{\a(\gamma-1)}{K_c(\gamma-3)(1+t)}
\phi^{\frac{\gamma-3}{2(\gamma-1)}}\right)(1+t)^{\frac{\alpha(3\gamma-1)}{2(\gamma-3)}},
 \end{aligned}
\end{equation}
and
\begin{equation}\label{4.2}
\begin{aligned}
{q_1}:=\left(\phi^{\frac{\gamma+1}{2(\gamma-1)}}{B}-\frac{\a(\gamma-1)}{K_c(\gamma-3)(1+t)}
\phi^{\frac{\gamma-3}{2(\gamma-1)}}\right)(1+t)^{\frac{\alpha(3\gamma-1)}{2(\gamma-3)}}.
 \end{aligned}
\end{equation}
Then $y_1$ and $q_1$ satisfy the following Riccati equations:

\begin{lemma}\ \label{L4.1}
 When $\lambda=1$, for any $C^1$ solutions of (\ref{1})-(\ref{4}),  we have
\begin{equation}\label{4.3}
\begin{aligned}
{y_1}^\prime =b_0-b_2 {y_1}^2,
\end{aligned}
\end{equation}
\begin{equation}\label{4.4}
\begin{aligned}
{q_1}^\backprime =b_0-b_2 {q_1}^2,
\end{aligned}
\end{equation}
where
\begin{equation}\label{4.5}
\begin{aligned}
b_0=\frac{  \a(\gamma-1)(\gamma-3)-\a^2(\gamma-1)^2}{K_c(\gamma-3)^2(1+t)^{2}}
\phi^{\frac{\gamma-3}{2(\gamma-1)}}(1+t)^{\frac{\alpha(3\gamma-1)}{2(\gamma-3)}},
\end{aligned}
\end{equation}
\begin{equation}\label{4.6}
\begin{aligned}
b_2=\frac{K_c(\gamma+1)}{2(\gamma-1)}
\phi^{-\frac{\gamma-3}{2(\gamma-1)}}(1+t)^{-\frac{\alpha(3\gamma-1)}{2(\gamma-3)}}.
\end{aligned}
\end{equation}
\end{lemma}

\textbf{Proof}.
By the proof in Lemma \ref{L3.1}, we have
\begin{equation}\label{4.7}
\begin{aligned}
\tilde{y}^\prime
=&-\tilde{a}_2\tilde{y}^2-\tilde{a}_1\tilde{y}+\tilde{a}_0,
 \end{aligned}
\end{equation}
where $\tilde{y}$ is defined by (\ref{26.6}), and $\tilde{a}_2$, $\tilde{a}_1$, $\tilde{a}_0$ are defined in (\ref{26.8.1}) by letting $\lambda=1$.
Then  we do one more simplification by multiplying
\begin{equation}\label{4.8}
\begin{aligned}
{\mu_1}=(1+t)^{\frac{\alpha(3\gamma-1)}{2(\gamma-3)}}
 \end{aligned}
\end{equation}
on (\ref{4.7}). In fact, it is easy to check that
\begin{equation}\label{4.9}
\begin{aligned}
{\mu_1}^\prime=\tilde{a}_1{\mu_1}.
 \end{aligned}
\end{equation}
Then  we denote
\begin{equation}\label{4.10}
\begin{aligned}
{y_1}={\mu_1}\tilde{y}.
 \end{aligned}
\end{equation}
Hence (\ref{4.7}) changes to
\begin{equation}\label{4.11}
\begin{aligned}
{y_1}^\prime=b_0-b_2{y_1^2}.
 \end{aligned}
\end{equation}
where
\begin{equation}\label{4.12}
\begin{aligned}
b_0={\mu_1}\tilde{a}_0, \ \ b_2=\frac{\tilde{a}_2}{{\mu_1}}.
 \end{aligned}
\end{equation}

Similarly, we have
\begin{equation}\label{4.13}
\begin{aligned}
\tilde{q}^\backprime
=&-\tilde{a}_2\tilde{q}^2-\tilde{a}_1\tilde{q}+\tilde{a}_0.
 \end{aligned}
\end{equation}
Again, if we denote $q_1=\tilde{\mu}\tilde{q}$,
then (\ref{4.13}) changes to
\begin{equation}\label{4.14}
\begin{aligned}
q_1^\backprime=b_0-b_2q_1^2.
 \end{aligned}
\end{equation}
This completes the proof.\\

\subsection{Uniform lower bounds on density for $1<\gamma<3$, $\lambda=1$}

\ \ \ \
In this subsection, the estimates of density's lower bounds will be considered for adiabatic exponent $1<\gamma<3$ and $\lambda=1$.
From the definition of $b_0$ in (\ref{4.5}), we have $b_0\l 0$. From Lemma \ref{L4.1}, we have
\begin{lemma}\ \label{L4.2}
When $1<\gamma<3$ and $\lambda=1$,
then for any $C^1$ solutions of  (\ref{1})-(\ref{4}), we have a priori bounds
\begin{equation}\label{E27}
\begin{aligned}
y_1(x, t)\leq \max\{1, \sup_x(y(x,0))\} =: Y,
\end{aligned}
\end{equation}
and
\begin{equation}\label{E28}
\begin{aligned}
q_1(x,t)\leq  \max \{1, \sup_x(q(x,0))\} =: Q.
\end{aligned}
\end{equation}
\end{lemma}

\begin{lemma}\ \label{L4.3}
When $1<\gamma<3$ and $\lambda=1$, assume system (\ref{1})-(\ref{4}) has a $C^1$ solution. Then there exists a positive constant $T>0$, such that
\begin{equation}\label{4.15}
\begin{aligned}
\rho(x,t)\geq K_0t^{-\frac{4}{3-\gamma}}(1+t)^{\frac{-2\alpha(3\gamma-1)}{(3-\gamma)^2}},\ \ \text{for any} \ \ t>T,
\end{aligned}
\end{equation}
where $K_0=\left[\left(\frac{2\sqrt{K\gamma}}{\gamma-1}\right)^{-\frac{\gamma-3}{2(\gamma-1)}}
\frac{3-\gamma}{2(\gamma-1)}K_c(Y+Q)\right]^{-\frac{4}{3-\gamma}}.$
\end{lemma}

\textbf{Proof}. From (\ref{20.1}), (\ref{4.1}) and Lemma \ref{L4.2},  we have
\begin{equation}\label{4.16}
\begin{aligned}
{w}_t&=-c{w}_x-\frac{\alpha}{2(1+t)}({z}+{w})\\
&=-K_c\phi^{\frac{\gamma+1}{\gamma-1}}
\left( y_1(1+t)^{-\frac{\alpha(3\gamma-1)}{2(\gamma-3)}}\phi^{-\frac{\gamma+1}{2(\gamma-1)}}
+\frac{\alpha(\gamma-1)}{K_c(\gamma-3)(1+t)}
\phi^{\frac{-2}{\gamma-1}} \right)-\frac{\alpha}{2(1+t)}({w}+{z})\\
&=-K_c\phi^{\frac{\gamma+1}{2(\gamma-1)}}
(1+t)^{-\frac{\alpha(3\gamma-1)}{2(\gamma-3)}}y_1-\frac{\alpha(\gamma-1)}{(\gamma-3)(1+t)}\phi-\frac{\alpha}{2(1+t)}({w}+{z})\\
&\geq-K_c\phi^{\frac{\gamma+1}{2(\gamma-1)}}(1+t)^{-\frac{\alpha(3\gamma-1)}{2(\gamma-3)}}Y
-\frac{\alpha(\gamma-1)}{(\gamma-3)(1+t)}\phi-\frac{\alpha}{2(1+t)}({w}+{z}).
\end{aligned}
\end{equation}
 Similarly, from (\ref{20}), (\ref{4.2})  and Lemma \ref{L4.2},  we have
\begin{equation}\label{4.17}
\begin{aligned}
{z}_t&=c{z}_x-\frac{\alpha}{2(1+t)}({z}+{w})\\
&=K_c\phi^{\frac{\gamma+1}{\gamma-1}}
\left(q_1(1+t)^{-\frac{\alpha(3\gamma-1)}{2(\gamma-3)}}\phi^{-\frac{\gamma+1}{2(\gamma-1)}}+\frac{\alpha(\gamma-1)}{K_c(\gamma-3)(1+t)}
\phi^{\frac{-2}{\gamma-1}} \right)-\frac{\alpha}{2(1+t)}({w}+{z})\\
&=K_c\phi^{\frac{\gamma+1}{2(\gamma-1)}}(1+t)^{-\frac{\alpha(3\gamma-1)}{2(\gamma-3)}}q_1
+\frac{\alpha(\gamma-1)}{(\gamma-3)(1+t)}\phi-\frac{\alpha}{2(1+t)}({w}+{z})\\
&\leq K_c\phi^{\frac{\gamma+1}{2(\gamma-1)}}(1+t)^{-\frac{\alpha(3\gamma-1)}{2(\gamma-3)}}Q
+\frac{\alpha(\gamma-1)}{(\gamma-3)(1+t)}\phi-\frac{\alpha}{2(1+t)}({w}+{z}),
\end{aligned}
\end{equation}
where $Y$ and $Q$ are defined in Lemma \ref{L4.2}. Therefore, for $1<\gamma<3$ we have
\begin{equation}\label{4.18}
\begin{aligned}
\phi_t&=\frac{1}{2}({w}_t-{z}_t)\\
&\geq-\frac{1}{2}K_c\phi^{\frac{\gamma+1}{2(\gamma-1)}}(1+t)^{-\frac{\alpha(3\gamma-1)}{2(\gamma-3)}}(Y+Q)
-\frac{\alpha(\gamma-1)}{(\gamma-3)(1+t)}\phi\\
&\geq-\frac{1}{2}K_c\phi^{\frac{\gamma+1}{2(\gamma-1)}} (1+t)^{\frac{\alpha(3\gamma-1)}{2(3-\gamma)}}(Y+Q).
\end{aligned}
\end{equation}
 Dividing the above inequality by $\phi^{\frac{\gamma+1}{2(\gamma-1)}}$, then integrating both sides on $t$,   we obtain
\begin{equation}\label{4.19}
\begin{aligned}
&\ \frac{2(\gamma-1)}{\gamma-3}(\phi(x,t))^{\frac{\gamma-3}{2(\gamma-1)}}
-\frac{2(\gamma-1)}{\gamma-3}(\phi(x,0))^{\frac{\gamma-3}{2(\gamma-1)}}\\
\geq&-\frac{1}{2}K_c(Y+Q)\int_0^t (1+s)^{\frac{\alpha(3\gamma-1)}{2(3-\gamma)}}ds\\
\geq&-\frac{1}{2}K_c(Y+Q)\int_0^t (1+t)^{\frac{\alpha(3\gamma-1)}{2(3-\gamma)}}ds\\
=&-\frac{1}{2}K_c(Y+Q)t (1+t)^{\frac{\alpha(3\gamma-1)}{2(3-\gamma)}}.
\end{aligned}
\end{equation}
Then
\begin{equation}\label{4.20}
\begin{aligned}
(\phi(x,t))^{\frac{\gamma-3}{2(\gamma-1)}}&\leq(\phi(x,0))^{\frac{\gamma-3}{2(\gamma-1)}}
-\frac{\gamma-3}{4(\gamma-1)}K_c(Y+Q)t (1+t)^{\frac{\alpha(3\gamma-1)}{2(3-\gamma)}}.
\end{aligned}
\end{equation}
There exists a positive constant $T$, when $t>T$,  we get
\begin{equation}\label{4.21}
\begin{aligned}
(\phi(x,t))^{\frac{\gamma-3}{2(\gamma-1)}}
\leq\frac{3-\gamma}{2(\gamma-1)}K_c(Y+Q)t (1+t)^{\frac{\alpha(3\gamma-1)}{2(3-\gamma)}},\ \ t>T.
\end{aligned}
\end{equation}
Therefore,
\begin{equation}\label{4.22}
\begin{aligned}
\rho(x,t)&\geq\left[\left(\frac{2\sqrt{K\gamma}}{\gamma-1}\right)^{-\frac{\gamma-3}{2(\gamma-1)}}
\frac{3-\gamma}{2(\gamma-1)}K_c(Y+Q)\right]^{-\frac{4}{3-\gamma}}
t^{-\frac{4}{3-\gamma}} (1+t)^{\frac{-2\alpha(3\gamma-1)}{(3-\gamma)^2}}\\
&:=K_0t^{-\frac{4}{3-\gamma}} (1+t)^{\frac{-2\alpha(3\gamma-1)}{(3-\gamma)^2}},\ \ t>T.
\end{aligned}
\end{equation}
This completes the proof.\\

\subsection{Singularity formation}

\ \ \ \
Similar to Section 3.3, this section is also divided into two cases \textbf{1)} $b_0\l 0$; and \textbf{2)} $b_0>0$. Moreover, in order to gain the
integration $\int_0^\i b_2(x(s),s)ds$,
the adiabatic exponent need to be separated into $\gamma>3$ and $1<\gamma<3$.\\

 \it \textit{\textbf{Case I:} For $\gamma>3$.}\\
\it From the definition of $b_0$ in (\ref{4.5}), We claim that $b_0\l 0$ is equivalent to
\begin{equation}\label{4.23}
 \a\g \frac{\gamma-3}{\gamma-1}.
\end{equation}
\it Using (\ref{4.6}) and Theorem \ref{t2.1}, we calculate
\begin{equation}\label{4.24}
\begin{aligned}
\int_0^\infty b_2(x(s),s)ds
&=\int_0^\infty\frac{K_c(\gamma+1)}{2(\gamma-1)}
\phi^{-\frac{\gamma-3}{2(\gamma-1)}}(1+s)^{-\frac{\alpha(3\gamma-1)}{2(\gamma-3)}}ds\\
&=\int_0^\infty\frac{K_c(\gamma+1)}{2(\gamma-1)}\left(\frac{2\sqrt{K\gamma}}{\gamma-1}\rho^{\frac{\gamma-1}{2}}\right)^
{-\frac{\gamma-3}{2(\gamma-1)}}(1+s)^{-\frac{\alpha(3\gamma-1)}{2(\gamma-3)}}ds\\
&=\frac{K_c(\gamma+1)}{2(\gamma-1)}\left(\frac{2\sqrt{K\gamma}}{\gamma-1}\right)^{-\frac{\gamma-3}{2(\gamma-1)}}
\int_0^\infty\rho^{\frac{3-\gamma}{4}}(1+s)^{-\frac{\alpha(3\gamma-1)}{2(\gamma-3)}}ds\\
&\geq\frac{K_c(\gamma+1)}{2(\gamma-1)}\left(\frac{2\sqrt{K\gamma}}{\gamma-1}\right)^{-\frac{\gamma-3}{2(\gamma-1)}}
\tilde{C}_0^{\frac{3-\gamma}{4}}\int_0^\infty (1+s)^{-\frac{\alpha(3\gamma-1)}{2(\gamma-3)}}ds,\\
&=K_1\frac{2(\gamma-3)}{2(\gamma-3)-\alpha(3\gamma-1)}(1+s)^{\frac{2(\gamma-3)-\alpha(3\gamma-1)}{2(\gamma-3)}}\big|^\infty_0\\
&=K_1\frac{2(\gamma-3)}{\alpha(3\gamma-1)-2(\gamma-3)}=K_1K_5>0,
\end{aligned}
\end{equation}
where $K_1=\disp\frac{K_c(\gamma+1)}{2(\gamma-1)}\left(\frac{2\sqrt{K\gamma}}{\gamma-1}\right)^{-\frac{\gamma-3}{2(\gamma-1)}}
\tilde{C}_0^{\frac{3-\gamma}{4}},
$ $K_5=\frac{2(\gamma-3)}{\alpha(3\gamma-1)-2(\gamma-3)}.$
In above calculation, we have used the function $2(\gamma-3)-\alpha(3\gamma-1)<0$.
 If we suppose the initial data
 $$y(x_0, 0)<-{1\over{K_1K_5}}:=-N_1,$$
 then $y(x(t), t)$ must blow up in finite time.\\

 Now, we give the singularity formation theorem for $\gamma>3$ in the following:
 \begin{theorem}\label{t4.1}(Singularity formation for $\gamma>3$)\\
\it Suppose the initial data $(\tau_0, u_0)(x)\in C^1(R)$ and there exists a positive constant $C_0$ such that
$$|u_0(x)|\l C_0,~ 0< \tau_0^{-1}\l C_0. $$
When $\a\g {{\gamma-3}\over{\gamma-1}}$,  assume there exists one point $x_0$ such that $y_1(x_0,0)\l -N_1$ or $q_1(x_0,0)\l -N_1$ , i.e.
\begin{equation}\label{4.28}
\begin{aligned}
&u_x(x_0,0)+\phi_x(x_0,0)<\tilde{K}_1(\phi(x_0,0))^{\frac{-2}{\gamma-1}}-N_1(\phi(x_0,0))^{-\frac{\gamma+1}{2(\gamma-1)}},
\end{aligned}
\end{equation}
or
\begin{equation}\label{4.29}
\begin{aligned}
&u_x(x_0,0)-\phi_x(x_0,0)<\tilde{K}_1(\phi(x_0,0))^{\frac{-2}{\gamma-1}}-N_1(\phi(x_0,0))^{-\frac{\gamma+1}{2(\gamma-1)}},
\end{aligned}
\end{equation}
where
$\tilde{K}_1=\frac{\a(\gamma-1)}{K_c(\gamma-3)},
~N_1={1\over{K_1K_5}},~\phi =\frac{2\sqrt{K\gamma}}{\gamma-1}\tau^{-\frac{\gamma-1}{2}}$. Then $u_x$ and/or $\tau_x$ must blow up in finite time.
\end{theorem}

\ \ \ \

\it\textit{\textbf{Case II:} For $1<\gamma<3$.}\\
 \it From Lemma \ref{4.2},  we have
\begin{equation}\label{4.50}
\begin{aligned}
(\rho(x,t))^{\frac{3-\gamma}{4}}\geq \left(K_0t^{-\frac{4}{3-\gamma}}(1+t)^{\frac{-2\alpha(3\gamma-1)}{(3-\gamma)^2}}\right)^{\frac{3-\gamma}{4}},\ \ t>T,
\end{aligned}
\end{equation}
then there exists a positive constant $T$ such that
\begin{equation}\label{4.51}
\begin{aligned}
&\ \ \ \ \int_0^\infty b_2(x(s),s)ds\\
&\geq\frac{K_c(\gamma+1)}{2(\gamma-1)}\left(\frac{2\sqrt{K\gamma}}{\gamma-1}\right)^{-\frac{\gamma-3}{2(\gamma-1)}}
\int_T^\infty\rho^{\frac{3-\gamma}{4}}(1+s)^{-\frac{\alpha(3\gamma-1)}{2(\gamma-3)}}ds,\\
&\geq\frac{K_c(\gamma+1)}{2(\gamma-1)}\left(\frac{2\sqrt{K\gamma}}{\gamma-1}\right)^{-\frac{\gamma-3}{2(\gamma-1)}}
\int_T^\infty\left(K_0s^{-\frac{4}{3-\gamma}} (1+s)^{\frac{-2\alpha(3\gamma-1)}{(3-\gamma)^2}}\right)
^{\frac{3-\gamma}{4}}(1+s)^{-\frac{\alpha(3\gamma-1)}{2(\gamma-3)}}ds\\
&=\frac{K_c(\gamma+1)}{2(\gamma-1)}\left(\frac{2\sqrt{K\gamma}}{\gamma-1}\right)^{-\frac{\gamma-3}{2(\gamma-1)}}K_0^{\frac{3-\gamma}{4}}
\int_T^\infty s^{-1}ds=\infty.
\end{aligned}
\end{equation}

\begin{theorem} (Singularity formation for $1<\gamma<3$)\label{t4.2}\\
\it Let the initial data $(\tau_0, u_0)(x)\in C^1(R)$ and there exists a positive constant $C_0$ such that
$$|u_0(x)|\l C_0,~ 0<\tau_0^{-1}\l C_0. $$
 Suppose  there exists one point $x_0$ such that $y(x_0,0)\l 0$ or $q(x_0,0)\l 0$ , i.e.
\begin{equation}\label{4.52}
\begin{aligned}
u_x(x_0,0)+\phi_x(x_0,0)<-\frac{\a(\gamma-1)}{K_c(3-\gamma)}(\phi(x_0,0))^{\frac{-2}{\gamma-1}},
\end{aligned}
\end{equation}
or
\begin{equation}\label{4.53}
\begin{aligned}
u_x(x_0,0)-\phi_x(x_0,0)<-\frac{\a(\gamma-1)}{K_c(3-\gamma)}(\phi(x_0,0))^{\frac{-2}{\gamma-1}},
\end{aligned}
\end{equation}
then $u_x$ and/or $\phi_x$ blow up in finite time.
\end{theorem}

\section{Comments and further problems}

\ \ \ \
To sum up this paper, we give some comments on our main results.\\
\it 1) From Theorem \ref{t3.1}-Theorem 4.2, we know that under the following two cases:\\
\it \it~~~~~~\textit{i) For $\gamma>3$, if $\ld \l \min\{1, {{\a(\gamma-1)}\over {\gamma-3}}\}$ or $\ld >\max\{1, {{\a(\gamma-1)}\over {\gamma-3}}\}$,}\\
\it \it~~~~~~\textit{ii) For $1<\gamma<3$ and $\ld \g{{\a(\gamma-1)}\over {\gamma-3}}$,}\\
the derivatives of Riemann invariants must blow up in finite time if their initial derivatives smaller than some fix constants, which only depend on the $L^\i$ norm of the initial Riemann invariants. In particular, $\ld<0$ means the damping coefficient increasing to infinity with the algebraic rate $\a(1+t)^{-\ld}$, the results in this paper show that the increasing damping can not cancel the nonlinear hyperbolic influence  totally.\\
\it 2) There have no the blow up results for $\ld\in \big(\min\big\{1, {{\a(\gamma-1)}\over {\gamma-3}}\big\}, \max\big\{1, {{\a(\gamma-1)}\over {\gamma-3}}\big\} \big)$ when $\gamma>3$ and $\ld<{{\a(\gamma-1)}\over {\gamma-3}}$ when $\gamma\in (1,3)$. However, we think this is just a technical problem, the derivatives of Riemann invariants may blow up in some certain initial data, too.\\
 \it 3) When   $ \lambda=0$, the system $(\ref{1})$ turns into the compressible Euler equations with constant coefficient damping.
 And our  results are also valid for this situation.

\section{Acknowledgement}
The authors would like to give many thanks to Professor Ming Mei and Ronghua Pan for their encouragement and discussion.

\end{document}